\documentclass{article}

\usepackage{epsf,latexsym,amssymb}
\usepackage{graphicx}

\newcommand{\edo}{\end{document}}

\newcommand{\R}{{\mathbb R}}  

\newcommand{\abs}[1]{\left\vert #1 \right\vert}

\newcommand{\be}[1]{\begin{equation}\label{#1}}
\newcommand{\ee}{\end{equation}}

\newcommand{\beq}{\begin{eqnarray}}
\newcommand{\eeq}{\end{eqnarray}}
\newcommand{\beqn}{\begin{eqnarray*}}
\newcommand{\eeqn}{\end{eqnarray*}}

\newcommand{\bi}{\begin{itemize}}
\newcommand{\ei}{\end{itemize}}

\newcommand{\ben}{\begin{enumerate}}
\newcommand{\een}{\end{enumerate}}

\newcommand{\pp}{p}

\title{\mbox{An observation regarding systems which converge to steady states}\\ for all constant inputs,
yet become chaotic with periodic inputs}
\author{Eduardo D.\ Sontag\\
Department of Mathematics\\
Rutgers University\\
New Brunswick, NJ 08903\\
{\tt sontag@math.rutgers.edu}}

\begin{document}

\maketitle

\begin{abstract}

This note provides a general construction, and gives a concrete example of,
forced ordinary differential equation systems that have these two properties:
(a) for each \emph{constant} input $u$, all solutions converge to a steady
state but (b) for some periodic inputs, the system has arbitrary (for example,
``chaotic'') behavior.  An alternative example has the property that
all solutions converge to the same state (independently of initial
conditions as well as input, so long as it is constant).

\end{abstract}

\section{The question}

We consider systems
\beqn
\dot x(t) = f(x(t),u(t))
\eeqn
of forced ordinary differential equations, where $u=u(t)$ is an input,
assumed scalar-valued for simplicity.

It is well-known fact that a system that exhibits periodic trajectories when
unforced may well exhibit chaotic behavior when forced by a periodic
signal, as shown by Van der Pol equations with nonlinear damping, and
nonlinear-force Duffing systems.  However, the following fact seems to be
less appreciated.

Assume that the system has the property that, whenever $u(t)$ is constant (any
value), all solutions converge to a steady state (which might depend on the
input and on the initial conditions).  Is it possible for such a system to
behave pathologically with respect to periodic inputs (no ``entrainment'')?
We outline a general construction together with a specific example of this
phenomenon.

A second construction is also given, which has the property that all solutions
for constant inputs converge to the same state (independently of initial
conditions as well as the input value, so long as it is constant).
This second construction displays a more interesting phenomenon, but
it is only done for an example (Lorentz attractor), and it seems difficult to
prove rigorously that its solutions when forced by a periodic input are
truly chaotic.

\section{The general construction}

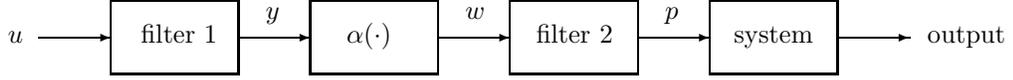
\begin{figure}[h,t]
\begin{center}
\setlength{\unitlength}{4144sp}%
\setlength{\unitlength}{2000sp}%
\begin{picture}(11280,924)(1336,-2998)
\put(1576,-2536){\vector( 1, 0){900}}
\put(2476,-2986){\framebox(1575,900){}}
\put(4051,-2536){\vector( 1, 0){900}}
\put(4951,-2986){\framebox(1575,900){}}
\put(6526,-2536){\vector( 1, 0){900}}
\put(7426,-2986){\framebox(1575,900){}}
\put(9001,-2536){\vector( 1, 0){900}}
\put(9901,-2986){\framebox(1575,900){}}
\put(11476,-2536){\vector( 1, 0){900}}
\put(5401,-2600){$\alpha (\cdot )$}%
\put(7750,-2600){filter 2}%
\put(10200,-2600){system}%
\put(12601,-2600){output}%
\put(1200,-2600){$u$}%
\put(2850,-2600){filter 1}%
\put(4396,-2311){$y$}%
\put(6871,-2311){$w$}%
\put(9346,-2311){$\pp$}%
\end{picture}%
\caption{General construction}
\label{fig:general}
\end{center}
\end{figure}

The general construction is illustrated in Figure~\ref{fig:general}.

The block called ``filter 1'' consists of a linear system $\dot x=Ax+Bu$ with
output $y=Cx+Du$ having the property that its transfer function $W(s)$ (the
Laplace transform of the impulse response $Ce^{tA}B$) is not identically zero
but has a zero at the origin: $W(0)=0$.

The differentiable function $\alpha :\R\rightarrow \R_{\geq 0}$ has the property that $\alpha (0)=0$
but $w(t) = \alpha (y(t))$ saturates at the value $1$ if $y(t)$ is outside a small
neighborhood of zero.

The block called ``filter 2'' consists of a scalar linear system
$\dot \pp = -\pp+w$.
If the signal $w(t)$ converges to zero, then $\pp(t)$ also converges to zero.
However, if the signal $w(t)\approx 1$ for a large set of times $t$, then
$p(t)$ will, after a transient, stay \emph{always} positive and bounded away
from zero, through an averaging effect.

Finally, the block named ``system'' has the form $\dot z(t) = \pp(t)f(z(t))$, where
$f$ is an arbitrary vector field with bounded solutions.  Provided that the
scalar signal $\pp(t)$ remains positive and bounded away from zero, the
solutions of the system will be time-rescaled versions of solutions of
$\dot z=f(z)$, and hence have the same pathologies.

When the input $u(t)$ is constant, the fact that the transfer function
vanishes at zero means that $y(t)\rightarrow 0$ as $t\rightarrow \infty $.
This implies that $w(t)\rightarrow 0$, and so $\pp(t)\rightarrow 0$ as well.
Thus the velocity $\dot z(t)\rightarrow 0$, and solutions become asymptotically constant.

On the other hand, suppose that the input is $u(t)=\sin \omega t$, where
$W(i\omega )\not= 0$.  This means that the output $y(t)$ will be asymptotically
of the form $A \sin(\omega t + \varphi)$, for some phase $\varphi\in \R$ and some amplitude
$A=\abs{W(i\omega )}\not= 0$.
Assuming that $A$ is large enough, passing $y(t)$ through $\alpha $ will result in
a signal $w$ that is mostly saturated near the value $1$, so $\pp(t)$ will stay
away from zero and the solutions of the $z$ subsystem will look (asymptotically)
like those of $\dot z=f(z)$.

\section{An explicit example}

Let us take ``filter 1'' to be the system $\dot x = -x-u$ with 
output $y=x+u$, the function $\alpha (y)=y^2/(K+y^2)$ with $K=1/10$,
and the vector field corresponding to the Lorentz attractor with
parameters in the chaotic regime.
Specifically, the complete system is as follows.
\def\lorx{\xi }
\def\lory{\psi }
\def\lorz{\zeta }
\beqn
\dot x&=&- x -u\\
\dot \pp&=&-\pp + \alpha (x+u)\\
\dot \lorx &=& \pp(s(\lory-\lorx))\\
\dot \lory &=& \pp(r \lorx- \lory- \lorx \lorz)\\
\dot \lorz &=& \pp(\lorx \lory - b \lorz)
\eeqn
with Prandtl number $s = 10$, Rayleigh number $r = 28$, and $b = 8/3$.
Figure~\ref{fig:example} shows typical solutions of this system with a
periodic and constant input respectively.

\begin{figure}[h,t]
\begin{center}
\includegraphics[scale=0.5]{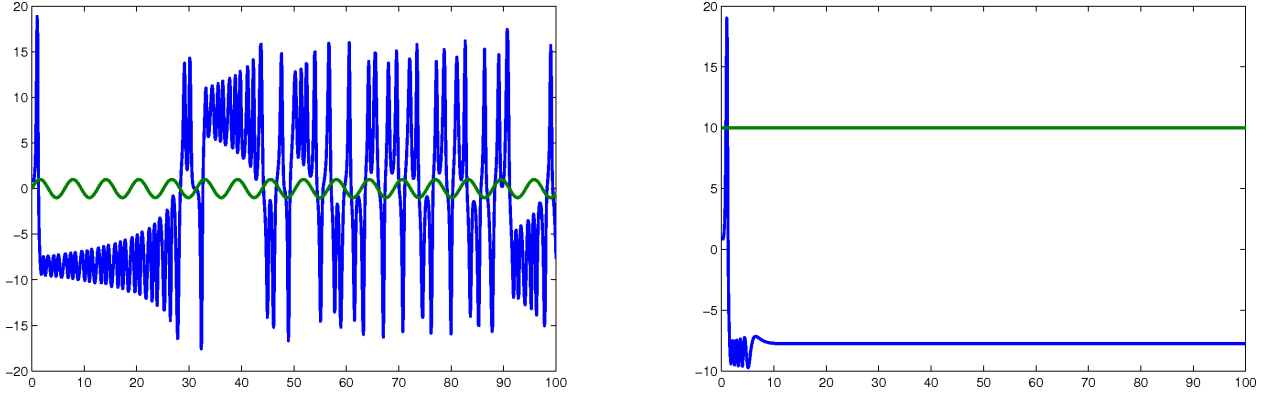}
\caption{Some solutions of example system.
Initial conditions are $x(0)=5$, $\lorx(0)=1$, and all remaining variables zero.
Green: inputs are $u(t)=\sin t$ (left panel) and $u(t)=10$ (right panel).
Blue: $\lorx(t)$.
Note chaotic-like behavior in response to periodic input, but steady state
in response to constant input.}
\label{fig:example}
\end{center}
\end{figure}

\section{An alternative example}

\begin{figure}[h,t]
\begin{center}
\includegraphics[scale=0.5]{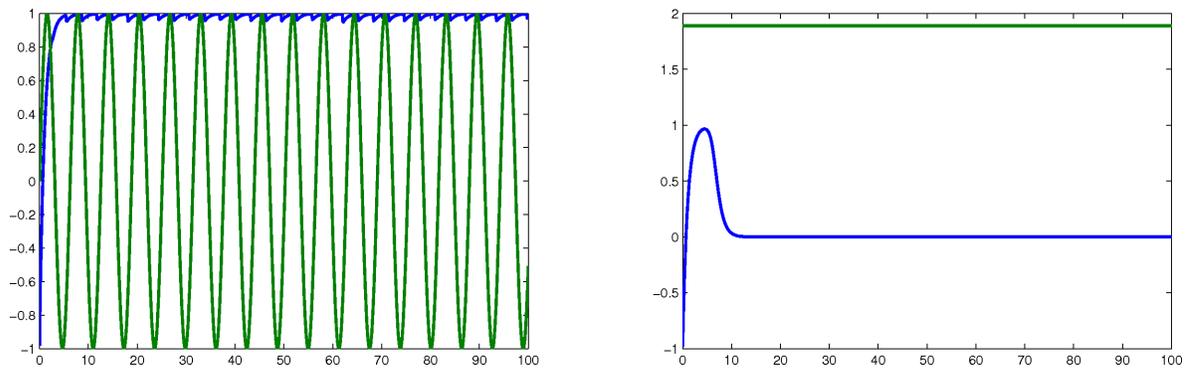}
\caption{Inputs and variable $\pp(t)$ for simulation shown
in~\protect{Figure~\ref{fig:example3}}.
Note that $\pp(t)\approx1$ and $\pp(t)\approx0$ for $t\gg0$.}
\label{fig:example3}
\end{center}
\end{figure}

Instead of multiplying the Lorentz vector field by $\pp(t)$, we may use $\pp(t)$
to eliminate the multiplicative terms, making the system asymptotically
linear and time-invariant, and also changing one of the coefficients in such a
way that the corresponding linear system is globally asymptotically stable.
Simulations (using a Montecarlo sampling of initial conditions as well as
input magnitudes) indicate that the system remains chaotic when forced by
$\sin t$, but all solutions converge to zero when $u(t)$ is constant.
Specifically, the complete system is as follows.
\def\lorx{\xi }
\def\lory{\psi }
\def\lorz{\zeta }
\beqn
\dot x&=&- x -u\\
\dot \pp&=&-\pp + \alpha (x+u)\\
\dot \lorx &=& 10(\lory-\lorx)\\
\dot \lory &=& 28\pp\lorx- \lory- \pp \lorx \lorz\\
\dot \lorz &=& \pp \lorx \lory - (8/3) \lorz
\eeqn
Observe that when $\pp=0$ one has, for the last three variables, the stable
linear system:
\beqn
\dot \lorx &=& 10(\lory-\lorx)\\
\dot \lory &=& - \lory\\
\dot \lorz &=& - (8/3) \lorz
\eeqn
and when $\pp=1$ the system is the usual Lorentz attractor
\beqn
\dot \lorx &=& 10(\lory-\lorx)\\
\dot \lory &=& 28\lorx- \lory- \lorx \lorz\\
\dot \lorz &=&  \lorx \lory - (8/3) \lorz \,.
\eeqn
In order to achieve $\pp\approx 1$ for the input $\sin t$, we now pick
$K=0.0001$ in the definition of $\alpha $.
Figure~\ref{fig:example2} shows typical solutions of this system with a
periodic and constant input respectively.
The function ``rand'' was used in MATLAB to produce random values in the
range $[-10,10]$.
\begin{figure}[h,t]
\begin{center}
\includegraphics[scale=0.5]{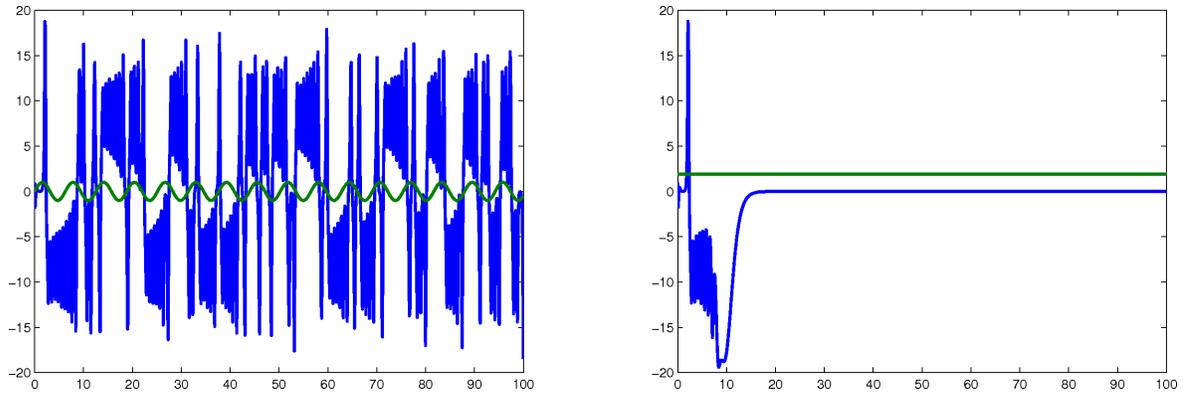}
\caption{Simulation of second example,
done with the following randomly-chosen input and initial conditions:
$u(t)=1.89$, 
$x(0)=2.95$ $\pp(0)=-0.98$, $\lorx(0)=0.94$, $\lory(0)=-4.07$, $\lorz(0)=4.89$.
Green: inputs are $u(t)=\sin t$ (left panel) and $u(t)=5.13$ (randomly picked,
right panel).
Blue: $\lorx(t)$.
Note chaotic-like behavior in response to periodic input, but steady state
in response to constant input.}
\label{fig:example2}
\end{center}
\end{figure}
To better understand the behavior of this system, Figure~\ref{fig:example3}
shows the value of the input and the variable $\pp(t)$, for the same
input and initial conditions as in Figure~\ref{fig:example2}.

\section{Remark}

One could clearly extend the construction to ``interpolate'' between any two
dynamics: using $\dot z(t)=p(t)f_1(z(t)) + (1-p(t))f_0(z(t))$ means that, for
constant inputs, asymptotically the system will behave like
$\dot z=f_0(z(t))$, while for $u(t)=\sin t$ it will behave like
$\dot z=f_1(z(t))$.  Establishing a precise result in this generality 
would require a detailed study of chain-recurrent sets and other
objects associated to asymptotically autonomous dynamics.
We leave that to future work.

\edo